

\baselineskip=14pt
\parskip=10pt

\magnification=\magstephalf

\def\1{{\overline{1}}}
\def\2{{\overline{2}}}

\parindent=0pt
\overfullrule=0in

\def\frac#1#2{{#1 \over #2}}
\centerline
{\bf  
Two Definite Integrals That Are Definitely (and Surprisingly!) Equal [Second Edition]
}
\medskip
\centerline
{\it Shalosh B. EKHAD, Doron ZEILBERGER and Wadim ZUDILIN}

\qquad \qquad \qquad \qquad {\it To our good friend Gert Almkvist (1934-2018), In Memoriam}

{\bf Proposition}: For  real $a>b>0$  and non-negative integer $n$, the following beautiful and surprising identity holds.
$$
\int_0^1 \, \frac{x^n\,(1-x)^n}{ ((x+a)(x+b))^{n+1}} \, dx \, = \,
\int_0^1 \, \frac{x^n\,(1-x)^n}{ ((a-b)x +(a+1)b)^{n+1}} \, dx \quad .
$$

{\bf Proof}: Fix $a$ and $b$, and let $L(n)$ and $R(n)$ be the {\it integrals} on the left and right sides respectively,
and let $F_1(n,x)$, and $F_2(n,x)$ be the corresponding {\it integrands}, so that
$L(n)=\int_{0}^{1} F_1(n,x)\, dx$ and  $R(n)=\int_{0}^{1} F_2(n,x)\, dx$.
We cleverly construct the rational functions
$$
R_1(x) \, = \,
{\frac {x \left( x-1 \right)  \left( (a+b+1) x^2+ 2abx -ab \right) }{ \left( x+b \right)  \left( x+a \right) }} \, , \,
R_2(x) \, = \,
{\frac {x \left( x-1 \right)  \left( (a-b)x^2+2b(a+1)x-(a+1)b \right) }{(a-b)x - (a+1)b}} \, ,
$$
with the motive that (check!)
$$
 \left( n+1 \right) F_1(n,x) - \left( 2\,n+3 \right)  \left( 2\,ba+a+b \right) F_1 \left( n+1,x \right) 
+ \left( a-b \right) ^{2} \left( n+2 \right) F_1 \left( n+2 ,x\right) \, = \,
\frac{d}{dx} (R_1(x)F_1(n,x)) \quad ,
$$
$$
 \left( n+1 \right) F_2(n,x) - \left( 2\,n+3 \right)  \left( 2\,ba+a+b \right) F_2 \left( n+1,x \right) 
+ \left( a-b \right) ^{2} \left( n+2 \right) F_2 \left( n+2 ,x\right) \, = \,
\frac{d}{dx} (R_2(x)F_2(n,x)) \quad .
$$
Integrating both identities from $x=0$ to $x=1$, and noting that the right sides vanish, we have
$$
 \left( n+1 \right) L \left( n \right) - \left( 2\,n+3 \right)  \left( 2\,ba+a+b \right) L \left( n+1 \right) 
+ \left( a-b \right) ^{2} \left( n+2 \right) L \left( n+2\right) \, = \, 0 \quad ,
$$
$$
 \left( n+1 \right) R \left( n \right) - \left( 2\,n+3 \right)  \left( 2\,ba+a+b \right) R \left( n+1 \right) 
+ \left( a-b \right) ^{2} \left( n+2 \right) R \left( n+2\right) \, = \, 0 \quad .
$$
Since $L(0)=R(0)$ and $L(1)=R(1)$ (check!), the proposition follows by mathematical induction. 

{\bf Comments}:{\bf 1}. This beatiful identity is equivalent to an identity buried in Bailey's classic book [B], section 9.5, formula (2), 
but you need an expert (like the third-named author) to realize that!
{\bf 2.} Our proof was obtained by the first named-author, running a Maple program, \hfill\break
{\tt http://sites.math.rutgers.edu/\~{}zeilberg/tokhniot/EKHAD.txt} , written by the second-named author, that
implements the Almkvist-Zeilberger algorithm [AZ] designed by Zeilberger and our good mutual friend Gert Almkvist,
to whose memory this note is dedicated. {\bf 3.} The integrals are not taken from a pool of no-one-cares analytic creatures: 
the right-hand side covers a famous sequence of rational approximations to $\log(1+(a-b)/((a+1)b))$ [AR]. Hence the left-hand side does.

{\bf Additional Comments (added in 2nd edition)}

{\bf 4.} This version corrects a sign typo pointed out by Greg Egan. {\bf 5.} Watch Greg Egan's beautiful animation in
{\tt https://twitter.com/gregeganSF/status/1192309179119104000}. {\bf 6.} To our surprise, the identity is not
as suprising as we believed. Mikael Sundquist noticed that the change of variable $x=b(1-u)/(b+u)$ gives
a `calc1 proof'. {\bf 7.} Alin Bostan has two further insightful proofs. See
{\tt https://specfun.inria.fr/bostan/publications/EZZ.pdf} and  {\tt https://specfun.inria.fr/bostan/publications/EZZ2.pdf}.

{\bf References}

[AR] K. Alladi and M. L. Robertson, {\it Legendre polynomials and irrationality},
J. Reine Angew. Math. {\bf 318}(1980), 137-155.

[AZ] Gert Almkvist and Doron Zeilberger,
{\it The method of differentiating under the
integral sign}, J. Symbolic Computation {\bf 10}(1990), 571-591; \hfill\break
{\tt http://www.math.rutgers.edu/\~{}zeilberg/mamarimY/duis.pdf}

[B] W. N. Bailey, {\it ``Generalized hypergeometric series''}, Cambridge University Press, 1935.

\bigskip
\hrule
\medskip
Shalosh B. Ekhad, Department of Mathematics, Rutgers University (New Brunswick), Hill Center-Busch Campus, 110 Frelinghuysen
Rd., Piscataway, NJ 08854-8019, USA. \hfill\break
Email: {\tt ShaloshBEkhad [at] gmail  [dot] com}   \quad .
\medskip
Doron Zeilberger, Department of Mathematics, Rutgers University (New Brunswick), Hill Center-Busch Campus, 110 Frelinghuysen
Rd., Piscataway, NJ 08854-8019, USA. \hfill\break
Email: {\tt DoronZeil [at] gmail  [dot] com}   \quad .
\medskip
Wadim Zudilin, Institute for Mathematics, Astrophysics and Particle Physics
Radboud Universiteit, PO Box 9010
6500 GL Nijmegen, The Netherlands
Email:  {\tt  wzudilin [at] gmail [dot] com} \quad .

\bigskip
\hrule
\bigskip
First Version: Nov.3, 2019. This version: Nov. 12, 2019.
\end